# From Plato's Rational Diameter to Proclus' Elegant Theorem

**by George Baloglou[1] and Yannis Thomaidis[2]**

**ABSTRACT:** We debate, departing from 2-3 readings of a single sentence in Proclus' Commentary to Plato's *Republic*, the plausibility of a rigorous (inductive) arithmetical derivation of an infinite sequence of pairs of side and diameter numbers by Proclus.

**ΠΕΡΙΛΗΨΗ:** Εξετάζουμε, ξεκινώντας από 2-3 αναγνώσεις μίας και μόνης πρότασης στα Σχόλια του Πρόκλου επί της *Πολιτείας* του Πλάτωνα, την πιθανότητα μιας αυστηρής (επαγωγικής) αριθμητικής παραγωγής μιας άπειρης ακολουθίας ζευγών πλευρικών και διαμετρικών αριθμών από τον Πρόκλο.

The Pythagoreans famously knew that there is no right isosceles triangle with integer sides: there are no integers a, d such that $d^2 = 2a^2$. The next best thing is to find integers a, d such that $d^2 = 2a^2 + 1$ or $d^2 = 2a^2 - 1$. Either the Pythagoreans or some of their epigones called integer pairs (a, d) with this property "pairs of side and diameter numbers" (***πλευρικοί και διαμετρικοί αριθμοί***), in obvious compensation for the fact that there is no right isosceles triangle with side length a and diameter (hypotenuse) length d. Indeed as these side and diameter numbers get bigger, the corresponding triangles get closer and closer to being right isosceles triangles:

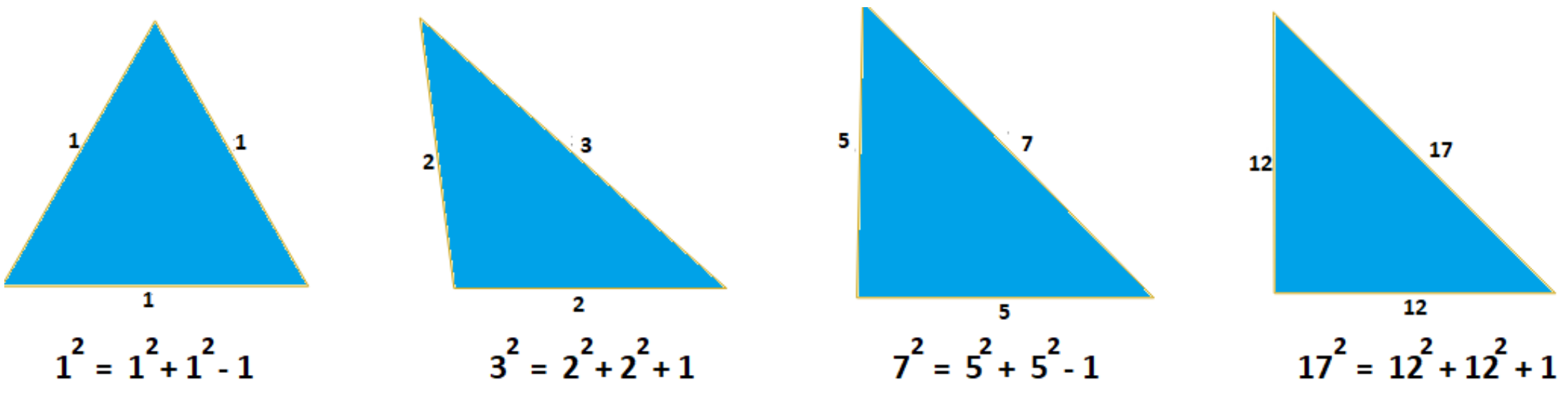

Figure 1

==================================================================

[1] Department of Mathematics, SUNY Oswego (1988 – 2008), USA – gbaloglou@gmail.com, PO BOX 50244, Thessaloniki 54006, GREECE

[2] Mathematics School Advisor (retired), Central Administration of Education – gthom54@gmail.com, Pesonton Iroon 1A, Thessaloniki 56430, GREECE

These four triangles were known to the Pythagoreans or at least Proclus and his predecessors, who also knew how to get from any such nearly right isosceles triangle to the 'next' one. The terms "algorithm" and "induction" naturally come to mind, even though such concepts were not known at the time. Still, Proclus' approach leaves room for speculation that he was interested in a 'general' and 'rigorous' way of getting from one pair of side and diameter numbers to the next and in a method of generating arbitrarily many such pairs; the possibility for such a goal is dramatically stressed by his concluding comment (end of section 27 in his Commentary to Plato's *Republic*) ***Και αεί ούτως*** ("And for ever like this"). Scholars differ regarding the level of 'rigor' sought and achieved by Proclus. Following an initially optimistic assessment [Heath 1921, van der Waerden 1954 & 1983, Knorr 1975], doubts began to emerge [Mueller 1981, Fowler 1987, Unguru 1991, Knorr 1998, Acerbi 2000]; optimism returned with [Negrepontis & Farmaki 2019].

In this paper we examine not only what Proclus does, but also what he *states* he does in a critical sentence (27.16-18 in his Commentary to Plato's *Republic* – [Kroll 1901, section 27 (KZ), p.27, 16-18]). As we are going to see, the optimistic reading of that key sentence relies on a mathematical leap involving the arithmetical reading of a geometrical proof and result, the so-called Elegant Theorem, allowed nonetheless by Proclus' own statements elsewhere. Specifically, Proclus views the arithmetical equality $d^2 = 2a^2 + 1$ or $d^2 = 2a^2 – 1$ as almost a substitute for the geometrical equality $d^2 = 2a^2$; further, this 'proximity' (***σύνεγγυς***) between the two equalities – that previous scholars did not relate to the issue under debate -- allows Proclus' geometrical argument associated with $d^2 = 2a^2$ to progress almost unchanged as an arithmetical argument associated with $d^2 = 2a^2 + 1$ and $d^2 = 2a^2 – 1$. In the optimistic reading the rigor of an ostensibly geometrical theorem and proof adds 'generality' to the arithmetical process that follows (specifically the first four examples of pairs of side and diameter numbers, each pair being derived from the previous one); in the established reading the said theorem serves mostly as a 'model' for the said process.

Note that the inductive and arithmetical reading(s) proposed in this paper are possible even under the interpretations of ***γραμμικώς*** and ***απ' εκείνου*** that are consistent with the traditional reading of Proclus' key sentence (27.16-18); these terms and issues are discussed extensively in Appendix II and in Appendices III & IV, respectively.

What could be Proclus' main motivation for 'strengthening', to one or another extent, in one or another way, the derivation of pairs of side and diameter numbers in a commentary to a mainly philosophical work? Probably not that work itself, but rather some otherwise unrecorded 'arithmetical maturity' of his era! The Elegant Theorem – viewable as a primitive mathematical induction step several centuries before al-Karaji [Rashed 1994], Gersonides [Rabinowitz 1970], and Maurolycus [Vacca 1909] – was known to other Neoplatonists, but Proclus was apparently the first one who felt the need for a rigorous proof related to it in one way or another: his ostensibly geometrical proof is flawless, a simple corollary of Euclid II.10; reversing the process in a way, Proclus also lauds the derivation from the Elegant Theorem of an arithmetical example of Euclid II.10, $3^2+7^2=2\times\left(2^2+5^2\right)$, apparently viewed by the Pythagoreans as a coveted substitute for the impossible $d^2 = 2a^2$.

## 1. An elegant theorem

"That's why Plato said about the number forty eight that it falls short of [the square of] the rational diameter of five by one [48 = 49 – 1] and of [the square of] the irrational [diameter of five] by two [48 = 50 – 2], because the [irrational] diameter's square is double of the side's square." [Kroll 1901, section 23 (ΚΓ), p.25, 6-9]

In a much debated passage (546c) in the *Republic*'s eighth book, Plato fatefully uses the term "rational diameter of five" for the number seven, on account of the 'nearly Pythagorean' equality $5^2+5^2=7^2+1$. This allowed for important testimonies on the early history of incommensurable quantities to reach us by way of Plato's commentators, notably Proclus of Lycia (Constantinople 412 CE – Athens 485 CE), also a noted Euclid commentator (Book 1 of the *Elements*).

According to Proclus [Kroll 1901, section 27 (ΚΖ), p.27, 1-6], the incommensurability between a square's side and diameter (diagonal) was from the beginning related to the impossibility of finding integer pairs (a, d) such that $d^2=2a^2$; after all, Plato also mentions, in the same *Republic* passage, the "irrational diameter of five". This impossibility led "the Pythagoreans and Plato" to invent "rational diameters" corresponding to "rational sides" [Kroll 1901, section 27 (ΚΖ), p.27, 6-10], that is pairs of integers (a, d) with squares satisfying the next best thing, namely either of the equations

$$d^2=2a^2-1 \quad \text{or} \quad d^2=2a^2+1.$$

Proclus gives [Kroll 1901, section 27 (ΚΖ), p.27, 10-11] concrete examples of such pairs (a, d) and their squares, that is (4, 9) and (25, 49). Further, he clearly suggests that the Pythagoreans were not limited to an empirical recording of a sequence of such pairs, but they devised recursive relations for producing successive terms of this sequence: if (a, d) is a term of the sequence, then $\left(a',\ d'\right)=\left(a+d,\ 2a+d\right)$ is another term of the sequence; in particular, from $3^2=2\times2^2+1$ one moves to $7^2=2\times5^2-1$, by way of 5 = 3 + 2 and 7 = 2x2 + 3. Indeed, right after the examples of squares (4, 9) and (25, 49), Proclus proclaims [Kroll 1901, section 27 (ΚΖ), p.27, 11-16]:

"And the Pythagoreans added (***προσετίθεσαν***) to this (***τούτου***) such an elegant theorem (***θεώρημα γλαφυρόν***) about the diameters and the sides, that the diameter increased by its side turns into a side [d + a = a'], and that the side added to itself and increased by its diameter turns into a diameter [a + a + d = d']."

As we are going to see, the choice of the geometrical terms "side" and "diameter" in relation to an arithmetical property is crucial, as the theorem just stated can be read both geometrically and arithmetically. This choice, by "the Pythagoreans and Plato" according to Proclus, was already made by the time of Theon of Smyrna (70 – 135 CE), who was also aware of the theorem, without stating it explicitly [Hiller 1878, pp.43-45]; likewise for Iamblichus (245 – 325 CE), who also introduces "side" and "diameter" in a clearly arithmetical context [Pistelli 1894, pp.91-92].

[We have corrected Kroll's ***προετίθεσαν*** ("placed before") to ***προσετίθεσαν*** ("added") … following the one and only surviving manuscript: see Appendix IV, Figure 4, line 7. On the other hand, we do preserve the standard rendering of ***γλαφυρόν*** as "elegant", already questioned in [Fowler 1987, 3.6 (b)], even though we argue in favor of "illuminating" or "effective" in Appendix I: such a meaning is more fitting with this ***θεώρημα γλαφυρόν*** being a tool rather than a goal, possibly a lemma needed for what would nowadays be called the "induction step", for moving from an *arbitrary* pair of side and diameter numbers to the next one.]

We note here that Proclus has already derived (25, 49) from (4, 9) – and (4, 9) from the "seed" (1, 1) – without explicitly stating the Elegant Theorem four sections earlier [Kroll 1901, section 23 (ΚΓ), p.24, 18-25 & p.25, 1-6]. But he does state, at the very beginning of section 23 (ΚΓ), that

"The Pythagoreans demonstrate arithmetically that the [squares of] rational [diameters] adjacent (***παρακείμεναι***) to irrational diameters are bigger or smaller than the double [of their sides' squares] by one". [Kroll 1901, section 23 (ΚΓ), p.24, 16-18]

For example, $\sqrt{8}$ is the irrational diameter of 2 whereas 3 is the adjacent rational diameter of 2, with $3^2 = 8 + 1 = 2 \times 2^2 + 1$; and $\sqrt{50}$ is the irrational diameter of 5 whereas 7 is the adjacent rational diameter of 5, with $7^2 = 50 - 1 = 2 \times 5^2 - 1$. (Recall here that the ancient Greeks had no way of *expressing* numbers like $\sqrt{8}$ or $\sqrt{50}$, hence the term ***άρρητοι αριθμοί*** ("non-expressible numbers") for irrational numbers, and likewise ***ρητοί αριθμοί*** ("expressible numbers") for rational numbers, both terms being related to ***ρήσις*** and/or ***ρήμα*** ("statement").)

Further down, Proclus writes:

“So if we wish to find a rational diameter, its square twice the square of five, taking the nearby (***σύνεγγυς***) with square not equal to it (for that is impossible), but smaller by one, we end up with 7, with a square of 49, which is smaller by one than 50, the square of the irrational diameter.” [Kroll 1901, section 35 (ΛΕ), p. 38, 17-22]

Back to the statement of the Elegant Theorem, 27.11-16 in [Kroll 1901, section 27 (KZ), p.27], what exactly does ***τούτου*** (“this”) stand for? It is unlikely to refer strictly to what precedes it (27.1-11), where the only item of direct relevance is the definition of “rational diameter” (***ρητή διάμετρος***) and the two examples (4, 9), (25, 49): Proclus has already stated a bit more in section 23, and in the statement cited above (24.16-18 in [Kroll 1901, section 23 (ΚΓ), p.24]) in particular; further, that ***Και αεί ούτως*** (“And for ever like this”) at the very end of section 27 is a clear indication about Proclus’ faith in the infinity of pairs of side and diameter numbers. We do not know to what extent the Pythagoreans – who, in Proclus’ words, “added to this such an elegant theorem” -- shared that faith; but it is reasonable to assume that anyone familiar with the Elegant Theorem would plausibly anticipate the said infinity. To put it differently, Proclus has very plausibly the ‘whole story’ (and in particular infinity) of pairs of side and diameter numbers in mind when he writes ***τούτου*** (“this”) in 27.12.

From here on, “property of rational diameters” stands for the statement that there exist infinitely many pairs of integers (a, d) such that $d^2 = 2a^2 - 1$ or $d^2 = 2a^2 + 1$; and Proclus’ ***τούτου*** (“this”) in [Kroll 1901, section 27 (KZ), p.27, 12] is assumed to refer to this statement.

## 2. What is “this”?

Right after the Elegant Theorem comes a short sentence [Kroll 1901, section 27 (KZ), p.27, 16-18] that we are going to interpret differently:

***Και τούτο δείκνυται δια των εν τω δευτέρω στοιχείων γραμμικώς απ’ εκείνου.***

“And this is demonstrated via the second [book]’s elements linearly from that.”

We explain in the next section and in Appendix II why we leave ***γραμμικώς*** essentially untranslated (“linearly”), hinting in favor of “rigorously” or “generally” instead of – or possibly in addition to -- previously established “geometrically”.

Proclus’ last two words above are ambiguous: ***εκείνου*** could be either masculine genitive of ***εκείνος*** = "he" or neutral genitive of ***εκείνο*** = "that"; and ***απ’*** = ***από*** could respectively mean either “by” or “from”.

Until very recently, historians bypassed this ambiguity in favor of "by him", "him" being none other than Euclid: under this interpretation, and since the Elegant Theorem follows indeed from Euclid II.10 (as Proclus shows right after the sentence in question and the statement of Euclid II.10), it must be credited to Euclid himself; in Fowler's rendering, "And this is demonstrated by lines (grammikos) through the things in the second [book] of Elements by him." [Fowler 1987, 3.6 (b)]

In this traditional "by him" interpretation there is the obvious issue of an uncalled for reference to someone who has not already been mentioned in the text, as well as the general lack of any attributions to Euclid – by other authors, that is – for results *not* included in the *Elements* (like the Elegant Theorem). These issues remain unresolved even after Heath's emendation [Heath 1926, p.400] of ***από*** (***απ'***) = "from" to ***υπό*** (***υπ'***) = "by", which most historians – Knorr in particular [Knorr 1975, p.33; Knorr 1998, p.422] – have tacitly followed since. Further, a TLG – [Thesaurus Linguae Graecae], a source used repeatedly here – search shows that ***από*** ("from") together with "prove" or "proof" ('roots' -***δεικ-***, -***δειξ-***, ***-δειχ-***) is used in Proclus' total corpus 29 times (always as "from", as for example in Proclus' Commentary to Euclid's First Book, ***το έβδομον δείκνυται από του πέμπτου*** = "Euclid I.7 is proven *from* Euclid I.5" [Friedlein 1873, p.269]), against only 5 times for ***υπό*** ("by"), in 4 of which Euclid is mentioned explicitly (***υπό του στοιχειωτού*** = "by the elements' creator"). (We also note here ***δείκνυται παρά τω Πτολεμαίω γραμμικώς*** = "is proven *by* Ptolemy linearly" [Manitius 1909, 4.106] as an alternative form of attribution employed by Proclus and other authors.)

[A closer look at [Heath 1926] throws more light into his emendation and interpretation. Heath actually misreads ***δια των εν τω δευτέρω στοιχείων*** ("via the second [book's] elements") as ***εν τω δευτέρω των στοιχείων*** ("in the second [book's] elements"), so he ends up reading Proclus as claiming that Euclid ("him") – and, by another leap, the Pythagoreans -- has proven the Elegant Theorem ("this") *in* his Second Book: "Proclus further says (p.27, 16-22) that the property of side- and diagonal- numbers "is proved graphically (***γραμμικώς***) in the second book of the Elements by him (***απ' εκείνου***)". Heath's probable motivation was to attribute the proof of the Elegant Theorem to the Pythagoreans (following [Zeuthen 1886]): "But the conjecture of Zeuthen, and the attribution of the whole theory of side- and diagonal- numbers to the Pythagoreans, have now been fully confirmed by the publication of "Procli Diadochi in Platonis rempublicam commentarii" (Teubner), Vol. II, 1901." [Heath 1926, p.398] (In textual reality, Proclus and others before him attribute the theory to the Pythagoreans, but Proclus is the only one who mentions (and then provides) a proof, with no specific attribution to anyone.)]

The preceding discussion establishes that we should opt for "from that" rather than "by him" at the end of the cited passage, choosing a fitting option, *seemingly* the only available, for "that" missed by previous researchers: we opt for ***εκείνου*** standing for the Elegant Theorem, that is the only result mentioned *so far* in section 27 (KZ) besides the property of rational diameters. This choice for ***εκείνου*** crucially implies in turn that "this" (***τούτο***) refers not to the previous sentence's Elegant Theorem, as historians assumed so far, but rather to the property of rational diameters preceding the Elegant Theorem [Kroll 1901, section 27 (KZ), p.27, 1-11] and the *first* "this" (***τούτου***): ***Προσετίθεσαν δε οι Πυθαγόρειοι τούτου τοιόνδε θεώρημα γλαφυρόν***……, "And to this [property of rational diameters] the Pythagoreans added such an elegant theorem……" [Kroll 1901, section 27 (KZ), p.27, 11-12].

This reading involves an obvious difficulty, "this" (***τούτο***), in 27.16-18 always, referring not to the immediately preceding Elegant Theorem (27.11-16) but to the property of rational diameters preceding the Elegant Theorem (27.1-11); likewise, "that" (***εκείνου***) referring to 27.11-16 rather than 27.1-11. It is perhaps this difficulty that motivated previous historians to associate ***τούτο*** with the Elegant Theorem and ***εκείνου***, *therefore*, with "him" (Euclid). But such usage of ***τούτο*** and ***εκείνου*** – and 'reversal' between them -- is far from impossible; see Appendix III.

It is tempting, perhaps, to associate ***εκείνου*** with ***δευτέρω***: "And this is demonstrated via the Elements' Second [Book] linearly from that [Second Book]". This reading has the 'proximity advantage' of associating ***τούτο*** with the preceding sentence's Elegant Theorem, but involves the obvious redundancy "demonstrated via the Second Book's elements from the Second Book".

More to the point, however, it has been pointed out to us that ***εκείνου*** may refer not necessarily to something preceding it, but possibly to something *following* it, in our case the very next sentence, which is none other than the statement of Euclid II.10! In this reading ***τούτο*** is naturally associated with the Elegant Theorem, while ***εκείνου*** comes as an explanation for ***δευτέρω*** (possibly with a slight (?) emendation involving a comma – ironically present in the manuscript after ***δείκνυται*** -- after ***γραμμικώς***): "And the Elegant Theorem is demonstrated via Euclid's Second Book, [that is] from Euclid II.10."; there is still an element of redundancy here, but not as strong as in the reading discussed in the previous paragraph. Such *cataphoric* use of ***εκείνου*** by Proclus – significantly (?) though involving a semicolon after ***εκείνο(υ)*** rather than a period as in 27.16-18 -- is discussed in Appendices III and IV.

[How 'natural' is the scheme "demonstrated via the second book, from that [theorem in the second book]"? We have no specific arguments against it, but would like to note here the apparent lack of combined use of "demonstrate" (***δείκνυται*** etc), "via" (***δια***) and "from" (***από***) in Greek.]

To summarize, we propose for further consideration the following two plausible readings for Proclus' key sentence (27.16-18):

(A) "And this [property of rational diameters] is demonstrated via the *Elements*' Second Book linearly from that [Elegant Theorem]."

(B) "And this [Elegant Theorem] is demonstrated via the *Elements*' Second Book linearly, from that [Euclid II.10]."

As we are going to discuss in the next section, reading (B) splits into two readings, depending on whether or not the Elegant Theorem – and the terms "side" and "diameter" -- is read arithmetically (reading (Ba)) or geometrically (reading (Bb)): since the arithmetical Elegant Theorem directly implies the property of rational diameters – inductively producing an infinity of pairs of side and diameter numbers – readings (A) and (Ba) effectively represent the same arithmetical reality; across them stands reading (Bb), representing a purely geometrical fact related nonetheless to the said arithmetical reality. (There is no point in similarly splitting reading (A), as the arithmetical Elegant Theorem yields the property of rational diameters at once; how and whether the geometrical Elegant Theorem could achieve the same goal is discussed in the next section.)

Of course (Bb) corresponds to the traditional reading, as presented for example in [Mueller 1981, p.288]: "… Proclus does not use II.9 and II.10 in this way. For him they are geometric results which can be proved exactly and to which the arithmetic procedure produces approximations. He makes no attempt to show that II.10 or the corresponding result about sides and diagonals of squares leads to an arithmetic truth when lines are interpreted as numbers. In other words, he does not take II.10 as an algebraic law."

Proclus' stance might have been a bit more complicated, as indicated by his thoughts on the theorems of Euclid's Second Book … in his Commentary to Euclid's First Book:

"Common to both sciences are the theorems regarding sections (such as Euclid presents in his second book), with the exception of cutting a segment in extreme and mean ratio [Euclid II.11]." [Morrow 1970, p.130; Friedlein 1873, p.60].

## 3. Geometrical or … rigorously arithmetical?

Proclus continues with a geometrical statement of Euclid II.10 [Kroll 1901, section 27 (KZ), p.27, 18-22] and a proof [Kroll 1901, section 27 (KZ), p.27, 22-24 & p.28, 1-9] of the Elegant Theorem's geometrical version from Euclid II.10.

[Fowler 1987, 3.6 (b)] records this material as follows:

“If a straight line be bisected and a straight line be added to it, the square on the whole line with the added straight line and the square on the latter by itself are together double the square on the half and of the square on the straight line made up of the half and the added straight line.”

[Euclid II.10: If A, B, Γ, Δ are collinear with AB = BΓ then $A\Delta^2 + \Gamma\Delta^2 = 2(AB^2 + B\Delta^2)$.]

“Let AB be a side and let BΓ be equal to it, and let ΓΔ be the diagonal of AB, having a square double that of it [i.e. AB]; by the theorem, the square on AΔ with that on ΔΓ, will be double that on AB and on BΔ. Of these, the square on ΔΓ is double that on AB; and it remains that the square on AΔ is double that on BΔ, for if as whole is to whole, so is what is taken away to what is taken away, the remainder, also, will be to the remainder as the whole is to the whole. Then the diagonal ΓΔ, receiving in addition the side BΓ, is a side; and AB, taking in addition itself, [i.e.] the [side] BΓ, and its own diagonal ΓΔ, is a diagonal; for it has a square double that of the side [ΔB].”

Euclid II.10 above is rendered as $(2a+d)^2 + d^2 = 2\left(a^2 + (a+d)^2\right)$ in modern algebraic terms (with 2a corresponding to the “whole line” AΓ and d to the “added line” ΓΔ), and the Elegant Theorem is rendered as $d^2 = 2a^2 \Rightarrow (2a+d)^2 = 2(a+d)^2$, in modern algebraic terms always.

According to Proclus’ own words ([Kroll 1901, section 27 (KZ), p.28, 4-6]), the Elegant Theorem follows from Euclid II.10 by way of Euclid V.19: if (u+v) : (y+x) = 2 and v : y = 2, then u : x = 2; here $u = (2a+d)^2$, $v = d^2 = 2a^2$, $x = (a+d)^2$, $y = a^2$. Proclus departs *not* from $(2a+d)^2 + d^2 = 2a^2 + 2(a+d)^2$ (applying first Euclid II.1 or Euclid II.2, that is distributivity) but from $(2a+d)^2 + d^2 = 2(a^2 + (a+d)^2)$, as the very statement of Euclid II.10 dictates. (“*The whole is twice the whole and what is taken away is twice what is taken away, so the remainder is twice the remainder*.”)

As it has been pointed out to us, it is tempting, for the contemporary reader at least, to apply distributivity first, departing thus from $(2a+d)^2 + d^2 = 2a^2 + 2(a+d)^2$, and from $d^2 = 2a^2$ (Elegant Theorem hypothesis) always, in order to arrive at $(2a+d)^2 = 2(a+d)^2$ (Elegant Theorem conclusion) … merely by applying Euclid’s Common Notion 3 (“*If equals are subtracted from equals, then the remainders are equal.”*). Fowler himself fell into this, criticizing Proclus’ appeal to Euclid V.19 as “unnecessarily and irrelevantly heavy-handed” [Fowler 1987, 3.6 (b)]; on the contrary though, the application of Euclid V.19 does not require any distributivity at all, so this approach might even have seemed simpler to Proclus and his contemporaries! (Fowler also mentions Euclid VII.11, a close relative of Euclid V.19 – in fact a special case of Euclid V.19 -- dealing with numbers rather than magnitudes, but that should be left out of discussion: indeed Proclus uses ***όλον προς όλον*** (neutral gender, as in ***μέγεθος*** =

“magnitude”) for “as the whole is to the whole”, exactly as in Euclid V.19, but Euclid VII.11 uses ***όλος προς όλον*** (masculine gender, as in ***αριθμός*** = “number”).)

Closing his lengthy geometrical parenthesis, Proclus concludes section 27 (KZ) with an “arithmetical proof” or at least “arithmetical demonstration” of the property of rational diameters (***δεικνύσθω δε επί των ρητών διαμέτρων αριθμητικώς***, “and let us [now] show about rational diameters arithmetically”), consisting in fact of a derivation of the first four instances of integer pairs (a, d) satisfying $d^2 = 2a^2 \pm 1$, followed by the statement ***Και αεί ούτως***, “And for ever like this” [Kroll 1901, section 27 (KZ), p.28, 10-27 & p.29, 1-4].

Rereading section 27, there are two obvious questions to ask, the first one already raised by Negrepontis & Farmaki 2019 (p.371):

(i) Why does Proclus resort to Euclid II.10 in order to prove a geometrical result that should be easy for him to prove more directly?

It is indeed easy to see – without appealing to Euclid II.10 or even to the Pythagorean Theorem for that matter – that from every right isosceles triangle of side a and diameter (hypotenuse) d we may obtain another right isosceles triangle of side a+d and diameter 2a+d; we illustrate this observation here by way of two ‘proofs without words’, the first of them at least 136 years old [Bergh 1886], very elementary and geometrical in nature (and probably ‘accessible’ to Proclus):

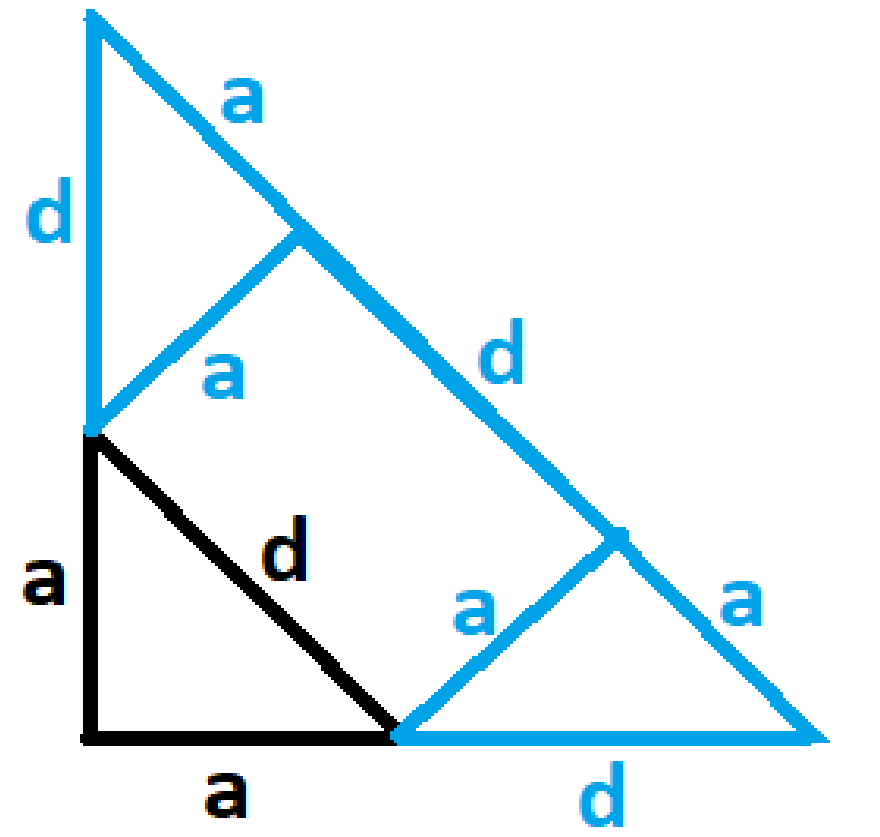

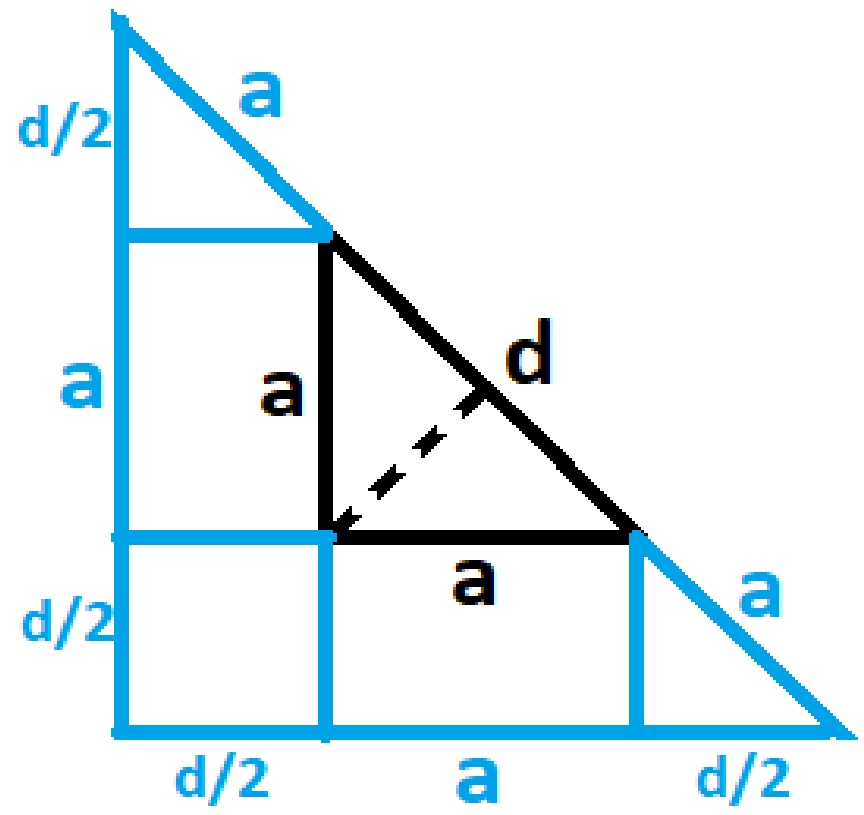

Figure 2

Assuming the Pythagorean Theorem, these visual proofs further establish the implication $d^2 = 2a^2 \Rightarrow (2a+d)^2 = 2(a+d)^2$ at once. Besides the possibility of Proclus (and others) simply not having ‘seen’ these proofs, there is also the plausibility of Proclus having deliberately chosen a ‘difficult’ proof *also* applicable to another, closely related, problem: this is discussed right below.

(ii) Why does Proclus cite a geometrical result and proof involving irrational diameters between two passages [Kroll 1901, section 27 (KZ), p.27, 1-11 and p.28, 10-27 & p.29, 1-4] devoted to arithmetical results about rational diameters? More to the point, are the Elegant Theorem [p.27, 11-16] and its proof [p.27, 22-24 & p.28, 1-9] formulated by Proclus as unequivocally geometrical entities?

Let us first observe that the first question in (ii) paves the way, together with the discussion in (i), towards a negative answer to the second question in (ii). Further, there is no reason to assume that the Pythagoreans, primarily interested in numbers rather than figures, would restrict themselves to an exclusively geometrical meaning of the terms "side" and "diameter" in the statement of the Elegant Theorem; and it is even less likely for Proclus himself, who begins section 27 with unequivocally arithmetical meaning for these terms, to limit himself to a strictly geometrical meaning only a few lines later.

Moreover, what the proof of the property of rational diameters and/or the proof of the Elegant Theorem read *arithmetically* (readings (A) and (Ba)) involves – in modern algebraic terms always -- is *not* the implication

$$d^2 = 2a^2 \Rightarrow (2a+d)^2 = 2(a+d)^2 ,$$

proven geometrically from Euclid II.10 for irrational diameters (Elegant Theorem read *geometrically*) by Proclus (reading (Bb)), but the 'adjacent' implication(s)

$$d^2 = 2a^2 - 1 \Rightarrow (2a+d)^2 = 2(a+d)^2 + 1 \text{ and } d^2 = 2a^2 + 1 \Rightarrow (2a+d)^2 = 2(a+d)^2 - 1$$

for rational diameters. Observe here that the statements of both implications are identical in the geometrical language used by Proclus ("the diameter turns into a side and the side turns into a diameter", without clarifying whether the diameter is irrational or rational) and that the 'algorithm' involved is the same (a' = d + a, d' = 2a + d, regardless of whether a and d are magnitudes or numbers).

At this point Proclus' statement on the 'proximity' between $d^2 = 2a^2$ and $d^2 = 2a^2 \pm 1$ necessary for our readings (A) and (Ba) and our application of Euclid V.19 is crucial: "…and when we are content with approximation, as for example in numbers, not having a square double of another square as in geometry, we say [the square of] a number smaller by one than twice [the square of] another number to be *double* of it, like [the square of] seven to [the square of] five"; this declaration from Proclus' Commentary to Euclid [Friedlein 1873, p.61] is strengthened by his statement regarding the Pythagorean Theorem (Euclid I.47), where he writes that there are no right isosceles triangles of integer sides … unless one considers $2a^2 \pm 1$ to be *double* of $a^2$ (***ει μη λέγοι τις τον σύνεγγυς***, "unless one means the nearby [number]") [Friedlein 1873, p.427].

Indeed 'allowing' $2a^2 \pm 1$ to be *double* of $a^2$, as Proclus is apparently inclined to do at times, crucially allows his proof (as rendered in [Fowler 1987, 3.6(b)] at the beginning of this section) to remain intact *word by word*, simply reading "the square on AΔ is double that on BΔ" as "$A\Delta^2 = 2B\Delta^2 + 1$" or as "$A\Delta^2 = 2B\Delta^2 - 1$" when "the square on ΔΓ is double that on AB" actually means "$\Delta\Gamma^2 = 2AB^2 - 1$" or "$\Delta\Gamma^2 = 2AB^2 + 1$", respectively.

Resorting to algebraic notation, Euclid V.19 applies as before, with $u = (2a+d)^2$, $v = d^2 = 2a^2 + 1$, $x = (a+d)^2$, $y = a^2$: 'allowing' $v : y = 2$ 'allows' in turn $u : x = 2$, so Proclus' readers had to merely 'guess' that a surplus of one (as in $v = 2y + 1$) implies a deficit of one (as in $u = 2x - 1$); such a guess would be facilitated by the arithmetical examples presented both before and after Proclus' proof. We feel that much less of a leap is required here compared to the traditional 'reconstructive' derivation of the arithmetical Elegant Theorem, as for example in [Van der Waerden 1954, p.126].

Following reading (Bb) and traditionally interpreting ***γραμμικώς*** as "geometrically", it appears at first that there is no problem (save for certain philological difficulties pointed out in section 2); as already pointed out in [Negrepontis & Farmaki 2019, p.369], however, there is then the issue of Proclus stressing the geometrical nature of the derivation of a geometrical result (Elegant Theorem) from another geometrical result (Euclid II.10).

Ironically, "geometrically" would not be a problem for readings (A) or (Ba): in reading (A) an arithmetical result (property of rational diameters) is derived from another arithmetical result (Elegant Theorem read arithmetically) proven 'geometrically' (from *arithmetical* Euclid II.10); and in reading (Ba) an arithmetical result (Elegant Theorem read arithmetically) is proven 'geometrically' from *arithmetical* Euclid II.10. In other words, in this scenario Proclus simply wishes to stress the geometric *appearance* of the derivation of an arithmetical result from another arithmetical result!

In another direction, we suggest, in support of readings (A) or (Ba), that Proclus' use of "linearly" (***γραμμικώς***) probably reveals his awe in front of a proof that is in fact non-geometrical yet *rigorous*. Indeed as we show in Appendix II ***γραμμικώς*** could be used as "rigorously" rather than "geometrically" (in Greek in general and in Proclus' writings in particular). This uncertainty about the exact meaning of ***γραμμικώς*** reflects on our rendering it here ambiguously and literally as "linearly".

In further support of "rigorously" over "geometrically", let us point out here that Proclus juxtaposes ***αριθμητικώς*** ("arithmetically") with ***γεωμετρικώς*** ("geometrically") [Kroll 1901, section 35 (ΛΕ), p.36, 3-4 and p.40, 1-2] only a few pages after juxtaposing ***αριθμητικώς*** with ***γραμμικώς*** in section 27 (ΚΖ): if ***γραμμικώς*** meant "geometrically" and nothing else, why not use ***γεωμετρικώς*** instead of ***γραμμικώς*** in the debated ***δείκνυται γραμμικώς*** passage [Kroll 1901, section 27 (KZ), p.27, 16-18], too?

One could still argue, in favor of reading (Bb), that ***γραμμικώς*** has here the meaning of "rigorously *and* geometrically": a geometrical proof is by definition rigorous, of course, but Proclus may wish to stress rigor here, especially in comparison to the arithmetical demonstration on rational diameters that follows [Kroll 1901, section 27 (KZ), p.28, 10-27 & p.29, 1-4]. (Compare with Proclus' "linear proof" in his Commentary to Euclid discussed in Appendix II.)

Reading (Bb) humbly suggests the following: even though Proclus no longer promises a proof of the property of rational diameters in 27.16-18, the geometrical proof (and algorithm) of the Elegant Theorem from Euclid II.10 right afterwards could be 'adapted' either into a full arithmetical proof of the property of rational diameters (by the reader, as described here) or, more realistically, into the derivation of the arithmetical examples from one another that follows (as in Proclus' text).

Reading (Ba) may appear to be the most plausible, involving  an arithmetical reading of the Elegant Theorem and a straightforward – to the contemporary reader *at least* -- derivation of the arithmetical Elegant Theorem from arithmetical Euclid II.10 paralleling Proclus' derivation of the Elegant Theorem from Euclid II.10; but it still involves the philological difficulties pointed out in section 2. Reading (A) avoids these difficulties but it states that the property of rational diameters follows from the Elegant Theorem, whereas it actually follows from an arithmetical reading of its proof, an arithmetical reading that establishes the arithmetical Elegant Theorem (which in turn trivially implies the property of rational diameters).

We summarize the three readings and their consequences in the following diagram:

Proclus' section 27 (KZ) [Kroll 1901, pp.27-29]

PRD arithmetical examples [28.10 - 29.4]

**27.16-18**

(A) **ETg - - - - -> PRD**

(Ba) **II.10 ---------> ETa**

**II.10 ---------> ETa**

$d^2 = 2a^2 +- 1$ -----> $(d+2a)^2 = 2(a+d)^2 -+ 1$

(Bb) **II.10 ---------> ETg**

**II.10 ---------> ETg** $d^2 = 2a^2$ -----> $(d+2a)^2 = 2(a+d)^2$

**27.22 - 28.9**

PRD = Property of Rational Diameters [27.1-11]
ETg = Elegant Theorem read geometrically [27.11-16]
ETa = Elegant Theorem read arithmetically [27.11-16]
II.10 = Euclid's Proposition II.10 [27.18-22]

Figure 3

In a way, the three readings discussed, each one with its own weak points, are not that different from each other. In all three of them there is a geometrical theorem that may also be read arithmetically, with a geometrical proof that relies on 'quasi-arithmetical' results such as Euclid II.10 and Euclid V.19, a theorem that is closely related to an easily verifiable arithmetical result and algorithm. The very presence of such an 'elegant' geometrical theorem amidst an unquestionably arithmetical discussion strongly suggests an auxiliary role for this theorem and its proof, supporting either an earthly inductive process (reading (Bb)) that could easily be described without it or a path-breaking inductive proof (readings (A) and (Ba)) that would be hard to present without a clumsy yet rigorous geometrical cover; in the first case we merely know how to get from a *specific* pair of side and diameter numbers to the next one, in the second case we have a *proof* that the process involved works for an *arbitrary* pair.

Proclus himself might be surprised by this debate … in view of what he writes [Friedlein 1873, p.60] in his Commentary to Euclid (right before his already cited statement about the 'proximity' between rational diameters and irrational diameters [Friedlein 1873, p. 61]), already quoted at the end of our Section 2; indeed Proclus is very comfortable with the notion of ostensibly geometrical propositions, such as Euclid II.1 – Euclid II.10 and Euclid II.12 & Euclid II.13, also accommodating an arithmetical meaning: in view of the proximity between the statements and proofs of the Elegant Theorem and the property of rational diameters discussed in this section, Proclus' statement strengthens our suspicion about the former serving as a cover for the latter. More to the point, it is indeed possible that Proclus did not distinguish between the Elegant Theorem's two forms, arithmetical and geometrical!

## 4. Demonstration versus proof

In a recent work that has been a major source of inspiration to us, Negrepontis & Farmaki [2019] traditionally read ***γραμμικώς*** as "geometrically", associate ***τούτο*** with the Elegant Theorem (without discussing ***εκείνου*** at all), and conclude (8.6.1, 8.6.5) that "the geometrical proposition Euclid II.10 is used geometrically to prove the Geometrical Elegant Theorem for the [irrational] diameter and arithmetically to prove the Arithmetical Elegant Theorem for the rational diameter ["property Pell" = property of rational diameters]". In support of their conclusion, they meticulously point out (8.6.2-8.6.4) strong language similarities – "recognized with a certain surprise confirming their expectation" (p.381) -- between Proclus' proof of the Elegant Theorem [Kroll 1901, section 27 (KZ), p.27, 22-24 & p.28, 1-9] and his derivation of the four pairs (1, 1), (2, 3), (5, 7), (7, 12) from one another [Kroll 1901, section 27 (KZ), p.28, 10-27 & p.29, 1-4]. It is indeed remarkable to see geometrical terms employed in an arithmetical demonstration, and fully consistent with the plausibility of a geometrical result illustrating an arithmetical proof promoted in this paper. (An even closer look at Proclus' text shows him to drift further and further apart from the Elegant Theorem's geometrical language as he progresses through his arithmetical examples.)

Where we still digress from [Negrepontis & Farmaki 2019] – and previous researchers as well -- is that, in our view, thoroughly discussed in the preceding two sections, Proclus is probably interested not in two theorems (both proven or at least provable from Euclid II.10) but *one*, merely using the Elegant Theorem and its proof as necessary 'illustration' for either a rigorous inductive proof of the property of rational diameters (readings (A) and (Ba)) or an inductive process for the derivation of arbitrarily many pairs of side and diameter numbers (reading (Bb)).

Interestingly, in the more recent [Negrepontis, Farmaki, and Brokou 202X] the authors appear to have come closer to our position, stating that "In the final step of synthesis, Proposition II.10, in an arithmetic form, immediately deduced from the geometrical, is used for the proof of every inductive step of the Pell equation, stated in a manner linguistically identical with the statement of the elegant theorem (so that the latter serves as the model for identical language for every inductive step in the proof of the Pell equation)": indeed the "Geometrical Elegant Theorem" is now only a means to an end, that is the "Arithmetical Elegant Theorem" (Pell Equation); and a geometrical proof is a model for an arithmetical demonstration or proof, whereas we view this geometrical proof as either an arithmetical proof in disguise further supported by the arithmetical demonstration that follows (readings (A) and (Ba)) or as a hint of rigor for the said arithmetical demonstration that parallels it (reading (Bb)).

Proclus devotes two sections of his *Republic* Commentary to the topic of rational diameters. At the beginning of section 23 (ΚΓ) he states that "the Pythagoreans demonstrate the property of rational diameters arithmetically" (***δια των αριθμών οι Πυθαγόρειοι δεικνύουσιν***), citing the first three pairs of side and diameter numbers; and he ends section 27 (ΚΖ) with an arithmetical demonstration by himself (***δεικνύσθω δε επί των ρητών διαμέτρων αριθμητικώς***) of the property of rational diameters, citing the first four pairs of side and diameter numbers. In both cases he shows in detail how to get from each pair to the next one, applying the algorithm suggested by the Elegant Theorem and/or Euclid II.10:

2 = 1 + 1 and 3 = 2x1 + 1, 5 = 2 + 3 and 7 = 2x2 + 3, 12 = 5 + 7 and 17 = 2x5 + 7.

Does Proclus, a seasoned Euclid commentator, feel that such a demonstration is on a par with the proofs in the *Elements*, that his ***Και αεί ούτως*** ("And for ever like this") – and his "quasi-general proof" ([Acerbi 2000]) -- is as convincing as Euclid's ***όπερ έδει δείξαι*** ("quod erat demonstrandum")? Probably not, and perhaps this is why he feels compelled to juxtapose his “arithmetical demonstration” – ***δε*** in ***δεικνύσθω δε επί των ρητών διαμέτρων αριθμητικώς*** – with his immediately preceding rigorous “linear demonstration” (***δείκνυται γραμμικώς***) of the Elegant Theorem (and plausibly the property of rational diameters); but possibly yes as well, citing perhaps the arithmetical demonstration right after the perceived rigorous proof as equivalent rather than merely illuminating. (One might argue here that Proclus stresses that he offers a rigorous proof for irrational diameters followed by a mere arithmetical demonstration for rational diameters, but we have already argued that his “linear demonstration” may actually refer to and work for *both* rational and irrational diameters; so we may well have a rigorous proof working for both rational and irrational diameters juxtaposed with a demonstration working only for rational diameters.)

Note at this point the citing of specific arithmetical examples together with rigorous geometrical proofs of the first ten theorems from the *Elements*’ Second Book in medieval texts (Arabic, Hebrew, and Latin) [Corry 2013]. Further, anonymous Greek commentators to Euclid’s Books II through XIII – effectively ‘continuators’ of Proclus – provide simple arithmetical examples without proofs of any kind for Euclid II.1-II.10 [Heiberg & Menge 1888, pp.227-248] and Euclid II.12 & II.13 [Heiberg & Menge 1888, pp.251-257], but a rigorous geometrical proof of Euclid II.11, which they show to be ‘non-arithmetical’, “non-demonstrable through numbers” (***ου δείκνυται δια ψήφων***) [Heiberg & Menge 1888, pp.248-251].

Note also the juxtaposition of ***γραμμικόν*** (“linear” as “general” or “rigorous”) with ***δια των αριθμών*** (“by numbers” as “by arithmetical examples”) by Pappus of Alexandria (290 – 350 CE) pointed out in Appendix II: that coexistence of arithmetical example with rigorous proof (even if by mere reference to the latter) strengthens the plausibility of Proclus’ arithmetical examples and “linear” proof referring to one and the same result, not the least because they are both called “demonstrations” (***δεικνύσθω αριθμητικώς***, ***δείκνυται γραμμικώς***) by Proclus!

It should be observed that, while Proclus explicitly attributes the arithmetical demonstration and the Elegant Theorem to the Pythagoreans, he writes nothing about the origins of the “linear” proof of the Elegant Theorem -- be it of the geometrical form (reading (Bb)) or the arithmetical form (readings (A) and (Ba)) -- that follows; all the ingredients were available to the Pythagoreans, according to Proclus at least, but there is certainly some distance – the culture of proof-by-example of those times mentioned above notwithstanding -- from the step-by-step arithmetical demonstration of a few pairs of side and diameter numbers to the rigorous proof of either the general ‘induction step’ (readings (A) and (Ba)) or the closely related geometrical result (reading (Bb)) that Proclus presents.

## 5. Euclid II.10

The four examples of pairs of side and diameter numbers displayed by Proclus, (1, 1), (2, 3), (5, 7) and (12, 17), were already known to Iamblichus (245 – 325 CE) and Theon of Smyrna (70 – 135 CE), who were also aware of the identity

$$d^2 + (2a + d)^2 = 2 \times \left( a^2 + (a + d)^2 \right),$$

describing it in about the same terms as Proclus [Kroll 1901, section 23 (ΚΓ), p.25, 9-11], as "all [two] diameters will be twice all [two] side*s* in power [of two]" [Pistelli 1894, p.93; Hiller 1878, p.44]. Iamblichus in particular writes

***ώστε αεί την διάμετρον δυνάμει διπλασίαν είναι της πλευράς, καθάπερ και επί των γραμμικών δείκνυται***,

"so that the diameter's square will always be twice the side's [square], exactly as it is demonstrated in Geometry" [Pistelli 1894, p.93].

We see a clear desire to view the identity $d^2 + (2a + d)^2 = 2 \times \left( a^2 + (a + d)^2 \right)$ as a certain substitute for the impossible for integer sides and diameters $d^2 = 2a^2$; further, the word ***αεί*** ("always") strongly indicates that Iamblichus was well aware of that identity's full generality, well beyond the identities

$$1^2 + 3^2 = 2 \times \left( 1^2 + 2^2 \right),\ 3^2 + 7^2 = 2 \times \left( 2^2 + 5^2 \right) \text{ and } 7^2 + 17^2 = 2 \times \left( 5^2 + 12^2 \right)$$

*implied* by the four examples above. On the one hand this is not terribly surprising, since the identity $d^2 + (2a + d)^2 = 2 \times \left( a^2 + (a + d)^2 \right)$ is nothing but an arithmetical version of Euclid II.10, and on the other hand it strongly suggests that Iamblichus, and very likely Theon as well, were aware of the Elegant Theorem and the algorithm a' = a + d, d' = 2a + d (already 'implicit' in Euclid II.10, which they do *not* mention), but never referred to it explicitly.

So, the property of rational diameters was well known among Neoplatonists, but it fell on Proclus to *plausibly* hint on a rigorous proof by way of Euclid II.10, illustrating it with the help of an 'elegant' geometrical theorem and its derivation from Euclid II.10. And it is very likely that such a rigorous arithmetical proof was not known to Theon and Iamblichus: they were certainly aware of the inductive process, they certainly knew how to get from one pair of side and diameter numbers to the next arithmetically, they knew how to apply the 'induction step', but they were probably *not* aware of any need to prove it (as Proclus *plausibly* does).

We note in particular that Iamblichus' ***καθάπερ*** ("exactly as") in ***καθάπερ και επί των γραμμικών δείκνυται*** concerns the result rather than the proof: Iamblichus is certainly aware of several instances of the equality $d^2+(2a+d)^2=2\times\left(a^2+(a+d)^2\right)$ and he is probably aware of how closely related that is to Euclid II.10; yet he is probably *not* aware of any 'general' derivation of

$$(d+2a)^2=2(a+d)^2\pm 1 \text{ from } d^2=2a^2\mp 1,$$

that leads, starting from $1^2=2\times 1^2-1$ and $3^2=2\times 2^2+1$, to $1^2+3^2=2\times\left(1^2+2^2\right)$ and other *special cases* of *arithmetical* Euclid II.10! It is of course quite ironic that, from this point of view, Euclid II.10 ends up being used in order to prove a few special cases of itself; more precisely, perhaps, we should say that geometrical Euclid II.10 ends up, by way of Proclus' proof of the Elegant Theorem, establishing a few special cases of arithmetical Euclid II.10.

Before Negrepontis & Farmaki [2019], researchers failed to notice the importance of the example $3^2+7^2=2\times\left(2^2+5^2\right)$ in the writings of Theon, Iamblichus, and Proclus. In [Knorr 1975], for example, Knorr writes

"then from a theorem like II. 10, it may be proved that

$$d_n^2+d_{n+1}^2=2\times\left(a_n^2+a_{n+1}^2\right)\text{”},$$

very shortly before mentioning the ratios $d_n : a_n$ as approximations of $\sqrt{2}$ (which none of the cited ancient authors mentions, besides a passing reference to "side and diameter ratio/(in)commensurability" (***πλευρικός και διαμετρικός λόγος***) by Theon [Hiller 1878, p.43] and Iamblichus [Pistelli 1894, p.91]), clearly giving the impression that this is his own reading (abandoned in [Knorr 1998]). Further, Fowler, unsure about the role of the Elegant Theorem, writes [Fowler 1987, 3.6 (b)] of an "inconsequential discussion of side and diameter numbers" in relation to Proclus' citing of "like nine together with 49 [is double] of 25 plus 4",

$3^2+7^2=2\times\left(2^2+5^2\right)$ (***οίον η εννέα μετά του μθ' της του κε' και δ'***),

and subsequent crucial testimony "this is why the Pythagoreans were encouraged by this method" (***διό και οι Πυθαγόρειοι εθάρρησαν τη μεθόδω***) [Kroll 1901, section 23 (ΚΓ), p.25, 12-13].

So, Proclus – but apparently not Theon or Iamblichus – credits the Pythagoreans with the discovery of the 'encouraging identity'

$$d^2+(2a+d)^2=2\times\left(a^2+(a+d)^2\right)$$

for *certain* a and d, which might even have been more important to them – as a substitute for the impossible $d^2 = 2a^2$ – than its *inductive* 'ingredient approximations'

$$d^2 = 2a^2 \mp 1 \text{ and } (2a+d)^2 = 2(a+d)^2 \pm 1 .$$

The said 'encouraging identity' is of course a special case of arithmetical Euclid II.10, and might even have led to the discovery of Euclid II.10 by the Pythagoreans. This view is supported by [Negrepontis & Farmaki 2021], where the authors further state (3.1) that "We conjecture that the Pythagorean attempt to turn the rational diameter into a real, geometric one marked the Birth of Induction". (Our readings (A) and (Ba) of 27.16-18 imply this conjecture and a proof – by Proclus rather than the Pythagoreans -- of what would nowadays be called "the general induction step", which in this case is none other than the Elegant Theorem.)

Knorr [1998] on the other hand points to very strong similarities between Euclid's 'unnecessarily complex' proof of II.10 and certain Mesopotamian diagrams; same about closely associated Euclid II.9, which he shows capable of proving the irrationality of $\sqrt{2}$, essentially rediscovering Stanley Tannenbaum's famous proof [Conway & Shipman 2013]. Van der Waerden [1954], Fowler [1987] and Negrepontis & Farmaki [2019] point to possible anthyphairetic origins of Euclid II.10, and everybody agrees that it is not used anywhere in the *Elements*. (Note however that Apollonius of Perga (262 - 190 BCE) uses Euclid II.9 and Euclid II.10 in the Third Book (sections 27 and 28) of *Conica* [Saito 1985/2004].)

Negrepontis & Farmaki [2019] ominously state (p.371), with the simplicity of Figure 2 (of our section 3) in mind, that "Euclid II.10 was *not* invented in order to prove the Elegant Theorem"; turning this on its head one may wonder whether Figure 2 actually led to the discovery of Euclid II.10 – 'simply' by adding $(2a+d)^2 = 2(a+d)^2$ to $d^2 = 2a^2$ – and possibly the property of rational diameters 25 centuries ago or so. More likely though, the Pythagoreans were not aware of Euclid II.10 (or a proof of the Elegant Theorem); they were thus 'encouraged' by the special cases of arithmetical Euclid II.10, obtained from the Elegant Theorem and yielding coveted substitutes for $d^2 = 2a^2$, but they were probably unaware of arithmetical Euclid II.10 in its full generality, much less of Euclid II.10 as geometrically stated in the *Elements*.

Before Proclus, Aristarchus of Samos (310 – 230 BCE) and Heron of Alexandria (10 – 70 CE) were aware of the approximations $\sqrt{2} \sim 7/5$ and $\sqrt{2} \sim 17/12$, respectively ([Heath 1913, p.378-79] and [Heiberg 1912, p.206]), indirectly provided by the Elegant Theorem and Euclid II.10. Heron's approximation was most likely achieved following the so-called "Babylonian Method", based on the well known algorithm

$$\left(\frac{d}{a}\right)' = \frac{1}{2}\left(\frac{d}{a} + \frac{2}{\frac{d}{a}}\right),$$

which he certainly was aware of [Heath 1921, vol. 2, pp.323-326]. Aristarchus' approximation could not possibly have been derived via this method (as the equation

$$\frac{7}{5} = \frac{1}{2}\left(x + \frac{2}{x}\right)$$

has no real solution), so it is very likely that it is indeed related to the Elegant Theorem and the induced slower algorithm

$$\left(\frac{d}{a}\right)' = \frac{2a+d}{a+d} = \frac{2+\frac{d}{a}}{1+\frac{d}{a}}.$$

(Both authors depart from a = 2, d = 3, in fact from a = 1, d = 1.)

Contrary to what is often assumed (see for example the otherwise informative and imaginative section on *side-and-diagonal-numbers* in [Hoyrup 2016]), there is no solid textual evidence that "the Pythagoreans and Plato" – or even Proclus for that matter -- were interested in rational approximations of $\sqrt{2}$ such as those of Aristarchus and Heron. It seems that what really mattered to them was the approximations

$$\text{of } d^2 = 2a^2 \quad \text{by} \quad d^2 = 2a^2 \pm 1$$

and, to Proclus at least, their infinity.

In the reverse direction, even though Heron did not need the Elegant Theorem to approximate $\sqrt{2}$, he was still interested in Euclid II.10 and the preceding nine theorems of Book II, sufficiently so to provide arithmetical – or at least "less geometrical" -- proofs of them, starting a lasting tradition that went beyond Greek Mathematics [Corry 2013]. And towards the other end of the Greek world, Barlaam of Calabria (1290 - 1348) has provided "arithmetical proofs of the [first ten] theorems proven geometrically in the *Elements*' Second [Book]" (***αριθμητική απόδειξις των γραμμικώς εν τω δευτέρω των στοιχείων αποδειχθέντων***) [Heiberg & Menge 1888, pp.725-738].

## Concluding Remarks

Proclus clearly shows how to inductively obtain a few pairs of side and diameter numbers employing the algorithm provided by the Elegant Theorem, which he attributes to the Pythagoreans. He also promises and then provides a rigorous ("linear") proof, of seemingly geometrical nature, of either the arithmetical or the geometrical form of the Elegant Theorem: in the first case we have a rigorous inductive derivation of an infinity of pairs of side and diameter numbers hidden behind an ostensibly geometrical proof; in the second case we have a geometrical theorem that supports a closely related arithmetical reality. It is not at all impossible that Proclus did not really distinguish between the two forms of the Elegant Theorem and their consequences; still, his use of "linearly" instead of "geometrically" does indicate a desire to reach a higher level of rigor and generality.

## Appendix I: "elegant"

The exact meaning of ***γλαφυρόν*** in Greek may be elusive, and certainly changing over time. Still used today as "vivid" or "illustrative", it meant "hollow(ed)" (from ***γλάφω*** = "carve") in Homer, applied to caves and harbors in particular; later on it acquired additional meanings such as "polished", "subtle", "refined", "elegant" [Liddell & Scott 1992, p.351]. Limiting our search to Proclus, we see that he has used it only four more times in his writings, three times in his Commentary to Euclid's First Book and once in his Commentary to Plato's *Parmenides*:

– In the first Euclidean instance [Friedlein 1873, p.72] Proclus writes that a proposition is "elementary" (***στοιχειώδες***) if it is both ***απλούν*** and ***χαρίεν*** and also has a wide range of applications yet is not important enough to the whole of science to be "fundamental" (***στοιχείον***), like the convergence of every triangle's heights, whereas a proposition that is neither related to several others nor exhibits anything ***γλαφυρόν*** and ***χαρίεν*** is not even "elementary". Interestingly, Morrow [1970, p.59] renders ***απλούν*** as "simple" and ***γλαφυρόν*** as "graceful", and ***χαρίεν*** as "elegant", leading to the redundancy "graceful and elegant" in "***γλαφυρόν*** and ***χαρίεν***" and losing the simplicity of ***απλούν***; we believe this happens because Morrow is biased in favor of ***γλαφυρόν*** = "graceful", and not flexible enough to adjust to the specifics of the given situation and opt for ***γλαφυρόν*** = "simple" or something along these lines.

– In the second Euclidean instance [Friedlein 1873, p.200] Proclus writes about "referring to the ancient writers" concerning the ***γλαφυρώτερα*** of the *Elements*' contents, "cutting down on their endless loquacity", as opposed to presenting the ***τεχνικώτερα*** and ***μεθόδων επιστημονικών εχόμενα*** "paying greater attention to the working out of fundamentals than to the variety of cases and lemmas": Proclus clearly juxtaposes the “simpler” or “more basic” (***γλαφυρώτερα***) contents with the “more technical ones” (***τεχνικώτερα***) and “relevant to scientific procedures” (***μεθόδων επιστημονικών εχόμενα***). Morrow [1970, p.157] sticks again to ***γλαφυρόν*** = “elegant”, which in this case seems to be out of place.

– In the third Euclidean instance [Friedlein 1873, pp.219-220] Proclus writes that it is ***γλαφυρόν*** that the equilateral triangle is constructed from one length, whereas the isosceles triangle is constructed from two lengths and the scalene triangle from three. Proclus chooses to stress this ***γλαφυρόν*** fact right after the respective ***πολυθρύλητα*** (“much spoken of”) constructions. Once again Morrow [1970, p.172] opts for “elegant”, but “clear” or “important” seems much more appropriate.

– In the Commentary to *Parmenides* Proclus suggests, comparing two different methods for seeking the truth, that it is more ***γλαφυρόν*** to do so by way of multiple hypotheses: ***ζητούντες τα πράγματα, μάλλον δια ταύτης ευρήσομεν ταληθές ή δι’ εκείνης, γλαφυρώτερον δια των πολλών τούτων υποθέσεων ανιχνεύοντες το ζητούμενον*** [Cousin 1864, p.1007]. Here the fitting meaning of ***γλαφυρόν*** is “accurate” [Morrow and Dillon 1987, p.357] or “effective”, with “elegant” not being appropriate at all.

This excursion through Proclus’ uses and respective meanings of ***γλαφυρόν*** (“simple”, “basic”, “clear”, “important”, “effective”, “accurate”) suggests a ‘compromise’ along the lines of “illuminating”: this is *not* an established meaning, so it would be helpful to locate examples along the lines of ***γλαφυρόν*** = “illuminating” in other authors. A rare yet great likely example occurs in a passage – that we will also encounter in Appendix II – from Plutarch, who mentions in *Marcellus* 14.9 the students of Eudoxus and Archytas as ***ποικίλλοντες τω γλαφυρώ γεωμετρίαν***, “making Geometry more effective” or “making Geometry more illuminating”, “solving problems not subject to logical and linear (rigorous) proofs by way of sensible and mechanical examples” [Ziegler 1968, 14.9]. (We disagree here – and so would Plato in view of section 14.11 – with B. Perrin’s well established yet rather awkward 1917 translation “making Geometry more subtle”, more precisely “they embellished Geometry with the subtleties of Mechanics” [Perseus Digital Library].)

Another possible example is Nicomachus Gerasinus’ reference to the ***γλαφυρά*** and ***ασφαλής*** generation of perfect numbers [Hoche 1866, Book 1, 16.4]: Nicomachus reproduces the well known formula for even perfect numbers (Euclid IX.36) … that one may view – in addition to being “secure” (***ασφαλής***) – as any of “simple”, “effective”, “illuminating”, “elegant”…

**Appendix II: "linearly"**

Our analysis below justifies our decision to essentially leave ***γραμμικῶς*** (from ***γραμμή*** = "line") untranslated in Proclus' key sentence ([Kroll 1901, section 27 (KZ), p.27, 16-18]), selecting the deliberately ambiguous and literal "linearly" over the previously established "geometrically", and discretely favoring "rigorously". (To be precise, "geometrically" first appears in [Knorr 1998]; before that we see the equivalent terms "graphically" [Heath 1921 and Thomas 1939], "by means of lines" [Knorr 1975], "by lines" [Fowler 1987, 3.6(b)], etc. Interestingly, van der Waerden [1954 & 1983] avoids mentioning the term altogether!)

Proclus uses "linear(ly)" a total of 20 times: twice in his Commentary to Plato's *Republic*, once in his Commentary to Plato's *Parmenides*, 5 times in his Commentary to Plato's *Timaeus*, 4 times in his Commentary to Euclid's *Elements*' First Book, and 8 times in *Hypotyposis Astronomicarum Positionum*, a work mainly discussing Ptolemy's *Almagest*.

In the Commentary to *Parmenides* "linearly" means "straight": "[souls] proceed forward in a line, but again turn back on themselves in a circle and return to their own first principles" [Morrow and Dillon 1987, p.471].

In all instances in *Timaeus* and two of four instances in Euclid, "linear" is in one way or another related to the "linear numbers": those are the positive integers exceeding 1, defined by Nicomachus as "starting at 2 and advancing always by adding the same unit distance" [Hoche 1866, 7.3, 8-11]; commenting on Nicomachus, Joannes Philoponus adds, about one century after Proclus, that "they imitate the line, having one dimension" [Giardina 1999, 29.2].

[The positive integers, "linear numbers", are therefore viewed as "sequential numbers", tempting one to render ***γραμμικῶς*** as "sequentially" in 27.16-18; that would fit particularly well with our reading (A), but is otherwise unattested.]

In the third Euclidean instance, "linear" means "one-dimensional": it is stated there that "a quantity is linear if and only if it can be partitioned by a point" [Friedlein 1873, p.121].

The fourth and last instance of "linear" in Proclus' Commentary to Euclid's First Book is the one that matters the most here: Proclus elaborates on the Fourth Euclidean Postulate, 'explaining' in two different ways why "all right angles are equal to one another" [Friedlein 1873, pp.188-191], and he presents his second explanation as a "linear proof". Does "linear" here mean "geometrical"? Proclus' second approach is certainly geometrical in nature, but why does he call it a "linear proof" instead of "geometrical proof"?

Both 'proofs' of the said equality of every two right angles are certainly geometrical in nature. The one called "linear" by Proclus is close in spirit to a typical Euclidean proof, consisting of a sequence of deductive steps, whereas the other one is kind of 'philosophical', supposedly following from first principles and definitions. The former is an example of "proof by inference" (***απόδειξις εκ τεκμηρίων επιχειρούσα***), whereas the latter is an example of "proof by definition" (***απόδειξις από των ορισμών μέσων το ζητούμενον δεικνύουσα***); this distinction is made clear by Proclus, who cites the proofs of Euclid I.32 and Euclid I.1 as examples of "proof by inference" (***απόδειξις από τεκμηρίου***) and "proof by cause" (***απόδειξις από αιτίας***), respectively [Friedlein 1873, p.206]. Further, Proclus writes [Friedlein 1873, p.69] about "all kinds of arguments", "some reaching validity by cause" (***τους μεν από των αιτίων λαμβάνοντας την πίστιν***) and "some departing from inference" (***τους δε από τεκμηρίων ωρμημένους***). We are in a position to suggest that in this difficult passage ***γραμμική απόδειξις*** means "rigorous proof".

Note that the meaning "linear" = "rigorous" is well established in Greek – from this term's first indisputable occurrence in the 1st century CE to Proclus' times and beyond. Indeed, in the *Marcellus* passage already discussed in Appendix I, Plutarch [Ziegler 1968, 14.9] juxtaposes "proofs by way of sensible and mechanical examples" (such as the one needed for the problem of doubling the cube) with "logical and linear proofs" (***λογικής και γραμμικής αποδείξεως***). Somewhat later, Galenus writes about "linear proofs" (***γραμμικάς αποδείξεις***) and "rhetorical arguments" (***ρητορικάς πίστεις***) [*De foetuum formatione libellus* , Kühn 1822, vol. 4, p.695] and "linear theorems" (***γραμμικοίς θεωρήμασιν***) he employs "following god's orders" [*De usu partium*, Kühn 1822, vol. 3, p.838], and calls Astronomy "the linear theory" (***η γραμμική θεωρία***), gradually developed from basic theorems [*De animi cuiuslibet peccatorum dignotione et curatione*, Kühn 1822, vol. 5, p.86]. Still later, Christian writers refer to "linear proofs" and St. Basil the Great writes about "linear and technical verbosity" (***γραμμικής και εντέχνου φλυαρίας***) [Giet 1968, Homily 3, section 3]. And the expression ***γραμμικαίς ανάγκαις*** = "via linear necessities" is used by Neoplatonists, from Porphyrius (3rd century CE) to Olympiodorus (6th century CE); the latter writes about Socrates showing methodically (***γραμμικαίς ανάγκαις***) how fairness brings happiness [Westerink 1970, 18.1].

Back to Proclus, we find 8 uses of "linearly" related to Ptolemy's work. It may appear at first that in these 8 cases "linearly" finally meant "geometrically", possibly related to applications of Geometry to Astronomy. This impression may be strengthened by the rarity of "geometrical" in Ptolemy's work: Ptolemy uses it in only two highly 'specialized' circumstances, ***γεωμετρικάς μεσότητας*** = "geometric means (averages)" [Jan 1895, section 25] and ***γεωμετρικόν*** = "topographical" [Grasshoff & Stückelberger 2006, Book 1, 2.2]. As it turns out, Ptolemy's "linear(ly)" means mostly

“trigonometrical(ly)”; see for example footnote 74 in [Sidoli 2020] for a discussion of closely related ***δια των γραμμών*** (“by lines”) and its usage in Ptolemy’s *Almagest*.

Even Ptolemy appears to have once used ***γραμμικώτερον*** (“more linear”) as “more rigorous”, as observed in section 2.2.1 (“***Δια των γραμμών*** and ***δια των αριθμών***”, “*By lines* and *by numbers*”) of [Sidoli 2004]: Ptolemy reproduces two Pythagorean proofs in Music Theory, calling the first one “more logical” and the second one “more linear” [Düring 1930, 1.5]; as Sidoli writes, “there is nothing geometrical about the second proof, which is a perfectly rigorous series of deductions” and “in this passage we should read ***γραμμικώτερον*** as meaning “more rigorous”” (p.112).

Further, even ***δια των γραμμών*** is probably used as “rigorously”, or at least “strictly geometrically”, by Ptolemy in a couple out of the 16 cases mentioned in [Sidoli 2020]: in [Heiberg 1898, p.42, 20], for example, Ptolemy writes – in relation to the difficulties of trisecting the angle of $1.5^0$ (in order to obtain the angle of $.5^0$) – that ***διά των γραμμών ου δίδοται πως*** (“cannot be found by geometrical methods”) [Toomer 1998, p.54]; and in [Toomer 1998, p.48], Ptolemy writes about “computing chord lengths by a strict geometrical method”, ***δια της εκ των γραμμών μεθοδικής αυτών συστάσεως*** [Heiberg 1898, p.32, 2-3]. As for ***γραμμικώς*** itself, in ***δείκνυσθαι γραμμικώς δια των υποκειμένων θεωρημάτων*** [Heiberg 1903, p.194, 2-3], “derived geometrically by means of the theorems [already] established” [Toomer 1998, p.411], preceded by “via lines only”, ***δια μόνων των γραμμών*** [Heiberg 1903, p.193, 19], “by purely geometrical methods” [Toomer 1998, p.410], there is room to speculate that “linearly” means more “rigorously” than “geometrically”.

Proclus uses “linear(ly)” as “geometrical(ly)”, more precisely “trigonometrical(ly)”, *only* when he writes about Ptolemy’s work. Conversely, Proclus uses ***γεωμετρικός*** 45 times in his Euclid Commentary, referring in particular to ***γεωμετρικαί δείξεις***, ***γεωμετρικαίς εφόδοις***, ***γεωμετρικών πίστεων***, all three meaning “geometrical proofs” or “geometrical methods”; and he also uses ***γεωμετρικός*** in the three Plato commentaries mentioned (9 times in *Republic*, 13 times in *Parmenides*, 70 times in *Timaeus*).

Likewise, Pappus of Alexandria (290 – 350 CE) uses ***γραμμικώς*** instead of ***γεωμετρικώς*** for “geometrically” *only* in his *Commentary to Ptolemy*, where ***γεωμετρικώς*** does not appear at all. Generally, non-astronomical uses of “linear(ly)” as “geometrical(ly)”, either before or after Proclus, are relatively uncommon; the first and only example before Proclus we find is Alexander of Aphrodisias’ proof (c. 300 CE) of a nebulous reflection property in his commentary on Aristotle’s *Meteorologica* [Hayduck 1899, p.144]. (At about the same time Iamblichus used the expressions ***επί των γραμμικών*** and ***εν γραμμικοίς*** (“in Geometry”, literally “on/in linears”), in his Commentary to Nicomachus Gerasinus’ *Introduction to Arithmetic* [Pistelli 1894, pp. 44, 59, 61, 92, 93].)

All these observations and statistics, together with some arguments in our section 3, provide considerable evidence in favor of ***γραμμικώς*** meaning "rigorously" rather than "geometrically" in Proclus' key sentence (27.16-18).

Are there are other possibilities for **γραμμικώς** beyond "rigorously" and "geometrically"?

As it turns out, this is a multifaceted word, with several alternate meanings:

Starting from Pappus, we see that he writes, in his *Collection* [Hultsch 1876, Book 4, p.270], of "linear problems", that is problems that involve curves ("lines") of degree higher than two – as opposed to "solid problems" (involving conic sections) and "planar problems" (involving lines and circles). A couple of centuries after Pappus, this usage survives in Asclepius of Tralles' commentary to Nicomachus' *Arithmetic*, where he writes about the difficulty of doubling the cube, resulting into appeal to linear and conical methods (***και οι μεν γραμμικώς εύρον, άλλοι δε κωνικώς***) [Tarán 1969, II.17]. (Notice the complete reversal of meaning of "linear(ly)" from Galenus to Asclepius; and how the latter uses ***γραμμικώς*** as "rigorously" in the same work, and in relation to the geometric mean and a special case of Euclid VII.19 [Tarán 1969, I.157].)

Next, both Galenus [*De usu partium*, Kühn 1822, vol. 4, p.20] and Cleomedes [Todd 1990, Book 1, section 5] use the expression ***μόνον ου γραμμικός*** = "all but linear", the latter referring to arguments in favor of the earth being spherical: this is more likely to mean "all but straightforward" than "all but rigorous", and not very likely to mean "all but geometrical"; sometimes it might indeed be hard to say whether "linear" means "straightforward" or "rigorous", like in the Asclepius example above!

More to the point, it is hard to tell whether "linear(ly)" means "general(ly)" or "rigorous(ly)" in Pappus' *Collection*'s Book 2 (that survives only in part, with Book 1 missing altogether); here "geometrical" is ruled out as Book 2 discusses rather obscure arithmetical problems involving large numbers and referring to work of Apollonius that has not survived. For example, Pappus writes in Proposition 18 "and the linear has been proven by Apollonius", ***το δε γραμμικόν υπό του Απολλωνίου δέδεικται*** [Hultsch 1876, Book 2, p.8, 27-28]; likewise, in Proposition 25 Pappus writes "this has been proven linearly by Apollonius", ***τούτο δε γραμμικώς Απολλώνιος απέδειξεν*** [Hultsch 1876, Book 2, p.18, 10-11]. In both examples, and in a few other instances in Book 2, ***γραμμικόν*** ("linear") is juxtaposed with ***φανερόν δια των αριθμών*** ("clear via numbers" = "clear via arithmetical example(s)") … in a way strongly reminiscent of Proclus' juxtaposition of ***γραμμικώς*** ("linearly") and ***αριθμητικώς*** ("arithmetically"). Remarkably, Pappus has even used ***δια των γραμμών*** ("by lines") in this spirit [Hultsch 1876, Book 2, p.4, 3-4], possibly 'inspired' by Ptolemy's usage mentioned above; see also [Sidoli 2004, pp.110-111].

Interestingly, Pappus also juxtaposes ***αριθμητικώς*** and ***δια των αριθμών*** with ***γραμμικώς*** on two occasions in his Commentary to Ptolemy's *Almagest* ([Rome 1931, p.45] and [Rome 1931, p.282], respectively).

Last but not least, in the *Republic* instance not discussed in this paper … "linear number" (***γραμμικόν αριθμόν***, [Kroll 1901, p.170, 23]) is none other than Plato's infamous "geometrical number": this is the one and only instance where Proclus uses "linear number" for Plato's 'nuptial number', as opposed to a total of 6 uses of "geometrical number" (3 times in his Commentary to *Republic* [Kroll 1901, p.36, 3; Kroll 1901, p.66, 7; Kroll 1901, p.70, 22], 2 times in his Commentary to *Timaeus* [Diehl 1904, p.275, 23; Diehl 1904, p.278, 3], once in his Commentary to Euclid [Friedlein 1873, p.23, 22]). Remarkably, no other Plato commentator discusses this elusive number! Plato himself mentions ***αριθμός γεωμετρικός*** ("geometrical number") only once, in *Republic* 546c (14 words after "rational diameter"): "And this entire geometrical number is determinative of this thing, of better and inferior births." [P. Shorey 1935, Perseus Digital Library]. See also [Baltzly, Miles, and Finamore 2022, pp.251-265].

**Appendix III: "that"**

As our discussion in section 2 clearly indicates, there is a philological issue in Proclus' key sentence ([Kroll 1901, section 27 (KZ), p.27, 16-18]) regarding ***απ' εκείνου*** ("from that") and the possibility (reading (B)) that it refers to next sentence's full statement of Euclid II.10: such *cataphoric* use of ***εκείνο(υ)*** ("that") is not unthinkable in Greek, even if rare in case the two sentences are separated by a period (as in Proclus' text) rather than a semicolon; we discuss this issue here (but, crucially, also in Appendix IV), focusing mostly on Proclus' usage (as transmitted to us by his scribes and editors), and on his most mathematical text (*Commentary to Euclid*, [Friedlein 1873]) in particular.

The closest to 27.16-18 example we find in Proclus' *Commentary to Euclid* occurs in 245.19-24: "But the hypothesis is composite in the common theorem about triangles and parallelograms with the same altitude. And both [hypothesis and conclusion] are composite as in this: "the diameters of circles and of ellipses bisect both the areas and the lines that contain the areas"." [Morrow 1970, p.289] The crucial part here is ***κατ' αμφότερα δε ως επ' εκείνου·*** ("and both [hypothesis and conclusion] as in that;"), with the theorem about circles and ellipses immediately following – unsurprisingly, and rather correctly, Morrow renders "that" as "this"!

The only other example in Proclus' *Commentary to Euclid* involving a semicolon occurs only a few lines earlier, in 245.4-7: "Interwoven theorems are such as can be divided into simple ones, like this [theorem]: “triangles and parallelograms with the same altitude are to one another as their bases”;” [Morrow 1970, p.289]. Again, Morrow renders ***οίον εκείνο το θεώρημα·*** as “like this [theorem]:” rather than “like that [theorem]:”.

There is no example involving a period, but there are two more examples involving a comma in Proclus’ *Commentary to Euclid*, both within the same sentence (326.15-22): “For to prove that the lines constructed within the triangle are shorter than the lines of the triangle outside them, *he requires the theorem* [***εκείνου δείται του θεωρήματος,***] ~~that~~ in any triangle two sides are greater than the third; and for showing that the angle they contain is greater than that contained by the outer lines, *he uses the proposition* [***εκείνο αυτώ συντελεί,***] ~~that~~ in any triangle the exterior angle is greater than the interior and opposite angle.” [Morrow 1970, pp.347-348]

Beyond Proclus’ *Commentary to Euclid*, are there any examples of cataphoric use of ***εκείνο(υ)*** where the two sentences are separated by a period, as in Proclus’ text under discussion? A first TLG search limited to “***εκείνο(υ).***” produced no such cataphoric examples in Proclus’ works (where, incidentally, “***εκείνο(υ).***” and “***εκείνο(υ)·***” appear 14 and 10 times in total, respectively). And such cataphoric use of “***εκείνο(υ).***” is generally rare in Greek, a bit less rare in late Greek (possibly influenced by Latin): we do find 3 cataphoric uses of “…***εκείνο.***” (“… that.”) in Claudius Aelianus’ *De natura animalium* (2.37.5, 9.54.10, 16.23.2); interestingly, the same work (by a Roman writing in Greek) contains 2 cataphoric uses of “…***εκείνο·***” (“… that;”), under rather analogous circumstances (5.1.18, 16.40.1) [Scholfield 1959].

Another, not entirely unrelated, issue raised in our section 2 concerns the ‘reversal’ between “this” (***τούτο(υ)***) and “that” (***εκείνο(υ)***) necessary for our reading (A) of 27.16-18; indeed this usage suggests some ‘unexpected proximity’ for “that”, also present when “that” refers cataphorically to something immediately following in the next sentence (as in our reading (B) of 27.16-18). We do not find examples of this reversal in Proclus’ *Commentary to Euclid* – a clear indication of scarcity – but there is an interesting example in Proclus’ Commentary to Plato’s *Timaeus* [Diehl 1903, vol. I, p.195, 12-18], which we discuss below starting with two translations separated by nearly 200 years:

“But as Socrates in the recapitulation of his polity asserts, that the cause of memory to us is the unusualness of the things which we hear, thus Critias, in what is here said, ascribes this cause to the age of children. For everything that occurs to children at first, appears to be unusual.” [Taylor 1820, pp.163-64]

“And just as Socrates in the summary of his constitution offers us as a reason for remembering the unfamiliarity of the things we hear, so Critias in this passage offers the youth of children. And it is likely that the former is also a reason for the latter case -- unfamiliarity a cause of children's memory. For all things seem unfamiliar when they are first encountered by children.” [Tarrant 2011, p.295]

It is remarkable that the critical half-sentence here, ***και έοικεν εκείνο και τούτου αίτιον είναι, η αήθεια της των παίδων μνήμης·*** (“and it is likely that the former is also a reason for the latter case -- unfamiliarity a cause of children's memory;“), is entirely missing from the first translation, and somewhat misrendered in the second translation (in the sense that “former” and “latter” need to be swapped, for what is a cause for memory is unfamiliarity, and in the preceding sentence, where the critical half-sentence refers, “unfamiliarity” comes *after* “remembering”).

The translators’ apparent unease with the cited passage is not accidental: “that” (***εκείνο(υ)***) is typically associated in English – and other contemporary languages, including Modern Greek -- with something more distant than what “this” (***τούτο(υ)***) refers to; but, as the Liddell & Scott lexicon clearly states for ***εκείνο***, “distant” may have to do not necessarily with position within the sentence, but with place in the mind, too:

> A. the person there, that person or thing: <u>generally with reference to what has gone immediately before</u>; but when ***οὗτος*** and ***ἐκεῖνος*** refer to two things before mentioned, ***ἐκεῖνος***, prop. belongs to the more remote, in time, place, or thought, ***οὗτος*** to the nearer. [Liddell & Scott 1992, p.505]

Indeed, in the passage from Proclus’ Commentary to *Timaeus* above, the main thing under discussion is memory, the unusualness – note here Taylor’s correction of ***αλήθεια*** (truth) to ***αήθεια*** – being just a cause for memory, therefore “more remote in thought”. And in Proclus’ passage under question, the main thing – and main topic of section 27 -- is the property of rational diameters, with the Elegant Theorem (on which it relies) being “more remote in thought”.

Another example is found in the work that contains the passage discussed in this paper, Proclus’ Commentary to Plato’s *Republic* [Kroll 1899, p.191, 18-25]:

“In these lines Socrates sets out very vividly that he separates poetry from that part of the mimetic class that produces illusions, and he says that it aims only at pleasure and the entertainment (*psychagegia*) of its audience. This is because the type of this *mimesis* that creates illusions is separated from the type that is representational, inasmuch as that [representational kind] considers the correctness of the imitation, but this one considers only the pleasure which the majority experience in illusion (*phantasia*).” [Baltzly, Miles, and Finamore 2018, p.302]

Here we see “that” being associated with the representational part (***εικαστικόν***) and “this” with the illusional part (***φανταστικόν***) in the second sentence’s second half, although the latter precedes the former in the second sentence’s first half (***και γαρ*** ***ταύτης της μιμήσεως το φανταστικόν απολείπεται του εικαστικού***): this happens because the main thing here is the illusional part, already under discussion by Socrates/Plato/Proclus [Kroll 1899, p.190, 26-27]; the previously discussed representational part [Kroll 1899, p.190, 2-25], is by now “more remote in thought”.

## Appendix IV: punctuation

Here we venture down to manuscript level for the sake of further exploring the significance of the period after ***απ' εκείνου***, already discussed in Appendix III. It could indeed be argued that the distance from semicolon to period in Proclus' *Commentary to Euclid*, 245.4-7 & 245.19-24 (***οίον εκείνο το θεώρημα·*** & ***κατ' αμφότερα δε ως επ' εκείνου·***), is ‘small’ and allows therefore for a cataphoric use of ***απ' εκείνου.*** in Proclus' key sentence [Kroll 1901, section 27 (KZ), p.27, 16-18]. As we shall see below, the answer to this question is not easy and … may well depend on the century!

Let us first take a look at the segment of the critical follie (42v) -- of the only available manuscript for Proclus' *Commentary to Plato's “Republic”*, Vatican 2197 (9th century), further discussed in [Baltzly, Miles, and Finamore 2022, pp.xi-xii], fully available online – that contains the passage of interest (27.6-18), from “the Pythagoreans and Plato” to “from that.” (and the beginning of Euclid II.10):

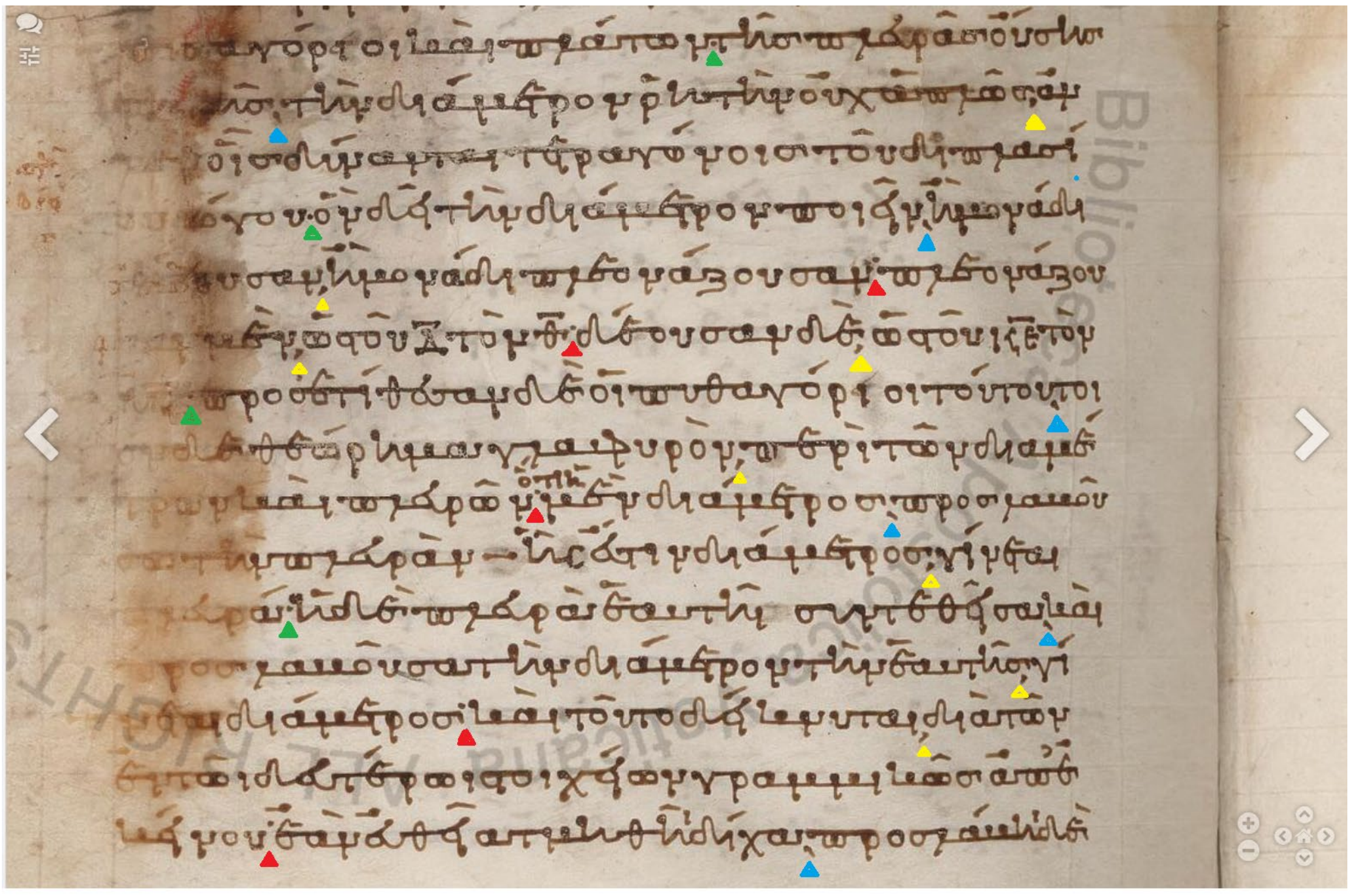

Figure 4

We indicate punctuation elements by triangles of various colors: yellow for commas, blue for 'reverse commas' (apparently corresponding to brief pauses), green for semicolon (looking like ano teleia), and red for period (looking like a further elevated ano teleia). So the period, not semicolon, right after ***απ' εκείνου*** (red triangle in the last line) is clearly visible and undeniable: Proclus – or at least the scribe(s) who produced this manuscript – could have had used a semicolon, but they opted for the period; still, there is a possibility that Proclus had used a semicolon that some scribe converted into a period!

Before moving further on, let us for comparison cite the text above as edited by Kroll in 1901:

μῆν τῶν γραμμῶν), ἐπενόησαν οὕτω λέγειν οἱ Πυθαγόρειοι
καὶ Πλάτων, τῆς πλευρᾶς οὔσης ῥητῆς τὴν διάμετρον ῥη-
τὴν οὐχ ἁπλῶς, ἀλλ᾽ ἐν οἷς δύνανται τετραγώνοις, τοῦ
διπλασίου λόγου, ὃν δεῖ τὴν διάμετρον ποιεῖν, ἢ μονάδι
10 [δέ]ουσαν ἢ μονάδι πλεονάζουσαν· πλεονάζουσαν μὲν ὡς
τοῦ Δ̄ τὸν Θ̄, δέουσαν δὲ ὡς τοῦ ΚΕ̄ τὸν ΜΘ̄. Προετί-
θεσαν δὲ οἱ Πυθαγόρειοι τούτου τοιόνδε θεώρημα γλαφυρὸν
περὶ τῶν διαμέτρων καὶ πλευρῶν, ὅτι ἡ μὲν διάμετρος
προσλαβοῦσα τὴν πλευράν, ἧς ἐστιν διάμετρος, γίνεται
15 πλευρά, ἡ δὲ πλευρὰ ἑαυτῇ συντεθεῖσα καὶ προσλαβοῦσα
τὴν διάμετρον τὴν ἑαυτῆς γίνεται διάμετρος. καὶ τοῦτο
δείκνυται διὰ τῶν ἐν τῷ δευτέρῳ στοιχείων γραμμικῶς
f. 43r. ἀπ᾽ ἐκείνου. ἐὰν εὐθεῖα τμηθῇ δίχα, προσλάβῃ δὲ | εὐ-

Figure 5

We also provide an English translation following the manuscript's punctuation:

"…… the Pythagoreans and Plato contrived the following wording; for a rational side … the diameter is [called] rational [to it] not as they are, but in their squares having a ratio of two; which must the diameter [in its square] produce … either missing a unit, or gaining a unit. Gaining, as 9 to 4. Missing, as 49 to 25; and the Pythagoreans added to this … such an elegant theorem, about the diameters and the sides. That the diameter … being added to the side to which is diameter, becomes a side; and the side added to itself … and increased by its own diameter, becomes a diameter. And this is demonstrated, via the second [book]'s elements linearly from that. If a line is cut in half … and is increased by ……"

Interestingly, even in this short fragment, there is a semicolon in the manuscript (green triangle in the first line) right after the cataphoric clause ***επενόησαν ούτω*** ***λέγειν Πυθαγόρειοι και Πλάτων*** ("the Pythagoreans and Plato contrived the following wording"), a semicolon that Kroll relaxed into a comma. Likewise, Kroll relaxed into a comma the manuscript's period (red triangle in line 9) following ***προσετίθεσαν δε οι Πυθαγόρειοι τούτου τοιόνδε θεώρημα γλαφυρόν περί των διαμέτρων και των πλευρών*** ("and the Pythagoreans added to this … such an elegant theorem, about the diameters and the sides"). But he did not opt for a similar 'relaxation' of the manuscript's period right after ***απ' εκείνου***: clearly, Kroll did *not* think of ***εκείνου*** as referring cataphorically to the following sentence (Euclid II.10); in the opposite direction, either Proclus or some scribe did *not* hesitate to 'separate' the Elegant Theorem's 'introduction' above from the theorem itself using a period!

This brief discussion about punctuation choices is quite inconclusive and by necessity very incomplete. A comprehensive study of punctuation in relation to cataphoric clauses would be a huge project way beyond the scope of this paper. Even if limited to Proclus, one should study all instances in all his works and all their manuscripts! Limiting ourselves to ***εκείνο(υ)*** in his *Commentary to Euclid* (as we did in Appendix III) and three of the best known witnesses to that work (one book and two manuscripts), we see – in Figures 6, 7 and 8 below, respectively – that Proclus (or at least some scribe) has mostly used commas or at least semicolons rather than periods in both 245.4-7 and 245.19-24:

Figure 6

[From page 66 in [Grynaeus 1533] (first printed edition).]

Figure 7

[From follies 99v & 100r of [Barocci 161 n.d.] (translated into Latin [Barocci 1560]).]

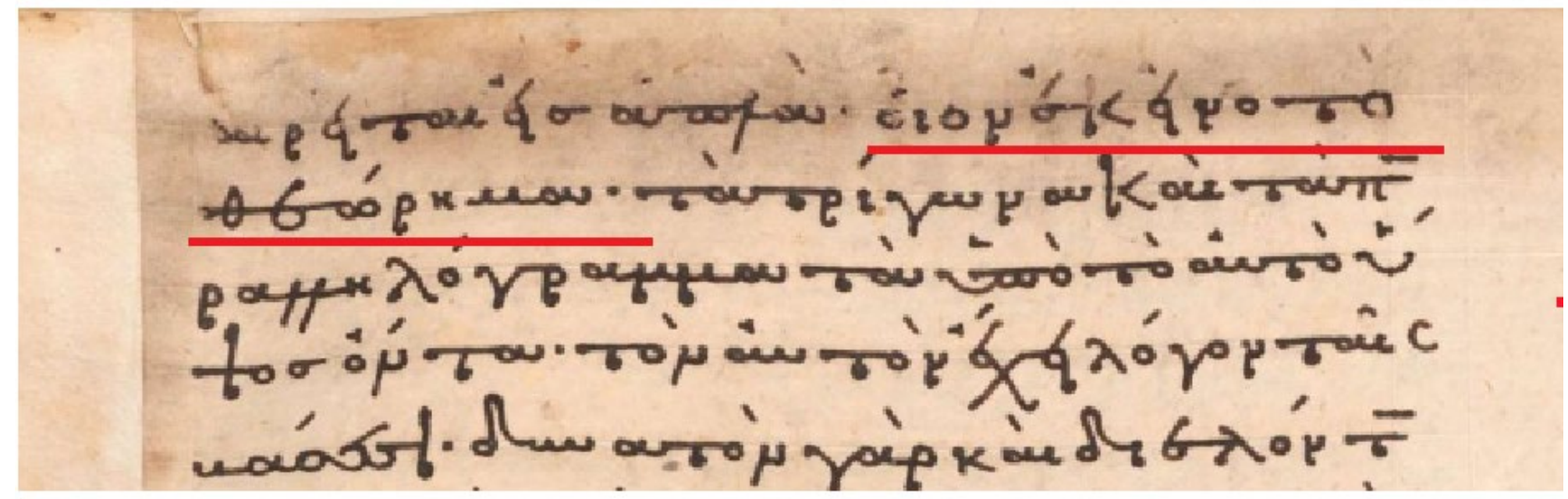

Fig. 8

[From follie 136v of [Monacensis 427 n.d.], main manuscript of [Friedlein 1873].]

In the first two works, that are relatively late (16th and 15th centuries), we see commas in both 245.4-7 and 245.19-24. In the third work (10th century) there is a period (again looking like a Modern Greek ano teleia, that is an ‘elevated period’) in 245.4-7 and a semicolon (looking like a period) in 245.19-24: so, even in a text where commas are rare anyway, and periods rather frequent, an anonymous scribe – or possibly Proclus himself -- opted for a semicolon rather than a period after ***ως επ’ εκείνου*** in 245.19-24, quite plausibly because of its clearly cataphoric use.

Clearly, semicolons in 245.4-7 and 245.19-24 of Proclus’ *Commentary to Euclid*, the only semicolons associated with cataphoric use of ***εκείνο(υ)*** in that work, are Friedlein’s choices, fully consistent with contemporary sensitivities, shared by Renaissance scribes and printers. (Incidentally, semicolons are totally absent from [Grynaeus 1533].) For the most part, cataphoric use of ***εκείνο(υ)*** appears to be associated with commas: this is also the case with Friedlein’s commas in 326.15-22 (discussed in Appendix III) and also in 102.16-18, 183.26-184.3, 203.23-204.2, again with the exception of the early manuscript. This of course does not necessarily imply non-cataphoric use when a period is used, as in our key sentence 27.16-18 of Proclus’ Commentary to Plato’s *Republic*: it is eminently possible that Proclus, or at least the 9th century scribe of [Vatican 2197], where/when periods are/were even more dominant than in [Monacensis 427], used a period in 27.16-18 even though they meant ***απ’ εκείνου*** to cataphorically refer to next sentence’s Euclid II.10 … and that [Kroll 1901] simply preserved that period assuming non-cataphorical use … plausibly misleading [Heath 1926] into associating ***απ’ εκείνου*** with Euclid (“him”)!

## Acknowledgements

We are grateful to Leo Corry, Theokritos Kouremenos, and Helen Perdicoyianni for their encouragement regarding an earlier version of this paper. We are also grateful to Alexander Alexakis, Alexandros Alexiou, Ioannis Avramidis, Dirk Baltzly, Dimitrios Christidis, Laval Hunsucker, Helen Karabatzaki, Stelios Negrepontis, Nick Nicholas, Katerina Oikonomopoulou, and Yannis Petrakis for helpful discussions and correspondence.